\documentclass[letterpaper, 10 pt, conference]{ieeeconf}  

\IEEEoverridecommandlockouts                              
\usepackage{graphicx}
\usepackage{amsmath, amssymb}
\usepackage{mathtools}
\usepackage[T1]{fontenc}
\usepackage{lmodern}
\usepackage{newtxtext, newtxmath}
\usepackage{subcaption}
\usepackage[noend]{algpseudocode}
\usepackage{algorithm}
\usepackage{nidanfloat}
\usepackage{url}
\usepackage{cite}
\usepackage{setspace}
\allowdisplaybreaks[1]

\newcommand{\T}{{\top}}

\newcommand{\pldev}[2]{{\partial {#1} / \partial {#2}}}
\newcommand{\bm}[1]{{\boldsymbol{\mathrm{#1}}}}
\newcommand{\R}{{\mathcal{R}}}

\newtheorem{theo}{Theorem}
\newtheorem{cor}{Corollary}
\newtheorem{defi}{Definition}
\newtheorem{assm}{Assumption}

\newtheorem{rem}{Remark}
\newtheorem{ex}{Example}

\algrenewcommand\algorithmicrequire{\textbf{Input:}}
\algrenewcommand\algorithmicensure{\textbf{Output:}}
\algrenewcommand\algorithmicforall{\textbf{for each}}

\title{\LARGE \bf
Bayesian Filtering for Nonlinear Stochastic Systems Using \\Holonomic Gradient Method with Integral Transform
}

\author{Tomoyuki Iori and Toshiyuki Ohtsuka
\thanks{This work was partly supported by JSPS KAKENHI Grant Numbers JP18J22093, JP21K21285, and JP15H02257.}
\thanks{T. Iori is with the Department of Information and Physical Sciences, Graduate School of Information Science and Technology, Osaka University, 1-5 Yamadaoka, Suita, Osaka 565--0871, Japan {\tt\small t-iori@ist.osaka-u.ac.jp}}%
\thanks{T. Ohtsuka is with the Department of Systems Science, Graduate School of Informatics, Kyoto University, Yoshida-Honmachi, Sakyo-ku, Kyoto 606--8501, Japan {\tt\small ohtsuka@i.kyoto-u.ac.jp}}%
\thanks{\copyright 2021 IEEE. Personal use of this material is permitted. Permission from IEEE must be obtained for all other uses, in any current or future media, including reprinting/republishing this material for advertising or promotional purposes, creating new collective works, for resale or redistribution to servers or lists, or reuse of any copyrighted component of this work in other works. }%
}

\begin{document}

\maketitle
\thispagestyle{empty}
\pagestyle{empty}

\begin{abstract}
   This paper proposes a symbolic-numeric Bayesian filtering method for a class of discrete-time nonlinear stochastic systems to achieve high accuracy with a relatively small online computational cost. 
   The proposed method is based on the holonomic gradient method (HGM), which is a symbolic-numeric method to evaluate integrals efficiently depending on several parameters. 
   By approximating the posterior probability density function (PDF) of the state as a Gaussian PDF, the update process of its mean and variance can be formulated as evaluations of several integrals that exactly take into account the nonlinearity of the system dynamics. 
   An integral transform is used to evaluate these integrals more efficiently using the HGM compared to our previous method. 
   Further, a numerical example is provided to demonstrate the efficiency of the proposed method over other existing methods. 
\end{abstract}

\section{Introduction \label{sec:intro}}
In most engineering applications, knowledge of the state of a dynamical system is required for monitoring or control via feedback. 
The optimal filtering theory was developed to estimate the current state of a system using the history of its outputs that are observable but contaminated by random noise. 
Following the success of the Kalman filter (KF)~\cite{Kalman1960b,Kalman1961} for linear systems with Gaussian noise, optimal filtering theory has been extended to applications involving nonlinear systems as well as non-Gaussian noise~\cite{Rusnak2015,Luo2015,Chen2018,Duong2018,Li2021}. 

Bayesian filtering is a general framework of optimal filtering. 
Although numerous problems can be formulated as the Bayesian filtering problems, they are difficult to implement as algorithms on computers, particularly for the real-time applications in systems and control. 
The Monte-Carlo scheme, which was first introduced for the particle filter (PF)~\cite{Kitagawa1993,Gordon1993} and subsequently used in other filtering algorithms, such as~\cite{Li2021}, is a powerful tool to compute integrations accompanied by statistical operations, such as marginalizations and expectations in nonlinear Bayesian filtering problems. 
Moreover, the nonlinearity of the system dynamics can be exactly accounted for by updating the particles according to the dynamics. 
To reap these benefits of the Monte-Carlo scheme, however, the number of particles must be sufficiently large, which requires heavy computational resources and may be unacceptable for applications with short sampling intervals. 
Another approach involves extending the KF to nonlinear cases such as the extended KF (EKF)~\cite{Jazwinski1970}, unscented KF (UKF)~\cite{Julier1997}, cubature KF (CKF)~\cite{Arasaratnam2009}, and Gauss-Hermite KF (GHKF)~\cite{Ito2000}. 
In all these methods, the posterior probability density function (PDF) is assumed to be or approximated as a Gaussian PDF. 
Under this Gaussian approximation, the marginalizations and expectations of the nonlinear functions, which are required to update the posterior mean and variance, can be efficiently computed using the Gauss quadrature rules. 
However, in exchange for efficiency, the Gauss quadrature rules cannot provide exact values of the integrals, which implies that the nonlinearity of the system dynamics can only be considered in the approximate sense. 
Thus, there is still room to extend the boundary of the trade-off between the exact consideration of nonlinearity and reduction of computational cost. 

In nonlinear cases, the integrals in the posterior PDF are ones of general nonlinear functions that depend on parameters such as the observed output, which are both analytically and numerically intractable to evaluate. 
In our previous work~\cite{Iori2020}, a symbolic-numeric method called the holonomic gradient method (HGM)~\cite{Nakayama2011} was used to evaluate the integrals efficiently.  
By utilizing the symbolic computation in terms of differential operators, we can perform integrations offline while considering the nonlinearity of the system exactly. 
The state estimation is then reduced to solving a set of initial-value problems (IVPs) of nonlinear ordinary differential equations (ODEs), which can be efficiently accomplished online with numerical integration methods such as the fourth order Runge-Kutta (RK4) or Adams-Bashforth-Moulton predictor-corrector (ABM4) methods. 
However, our previous method requires solving an IVP for each component of the posterior mean and variance in the one-step estimation. 
This incurs a high computational cost for the online computation and is unacceptable for short sampling intervals. 

To overcome the computational limitations of our previous method, we use an integral transform. 
Integral transforms are useful tools to compute the moments of a probability distribution efficiently. 
For example, in~\cite{Idan2014}, Idan and Speyer used an integral transform called the characteristic function of the multivariate Cauchy distribution. 
For linear systems under additive Cauchy noise, the characteristic function can be analytically propagated in time and can be used to derive the explicit forms of the posterior mean and variance. 
In this paper, we introduce an integral transform, which is different from, but similar to, the moment generating function. 
This similarity enables us to compute all the components of the posterior mean and variance efficiently from the integral transform and its derivatives. 
Moreover, the integral transform and its derivatives can be simultaneously evaluated by using the HGM only once. 
Hence, in contrast to our previous method, only one IVP must be solved, which leads to the reduction of the online computational cost. 

\paragraph*{Notations}
For the field of real numbers $\bm{R}$ and a vector of indeterminates $X = [X_1\ \cdots\ X_n]^\T$, $\bm{R}(X)$ denotes the field of rational functions in the components of $X$ over $\bm{R}$. 
$\partial_X \coloneqq [\partial_{X_1}\ \cdots\ \partial_{X_n}]^\T$ denotes a vector of differential operators, where $\partial_{X_i} = \pldev{}{X_i}$. 
We abbreviate $\partial_{X_i}$ by $\partial_i$ if $X$ is clearly specified according to the context. 
For a multi-index vector $d = [d_1\ \cdots\ d_n]^\T \in \bm{Z}^n_{\ge 0}$, $X^d$ and $\partial^d$ denote $X_1^{d_1} \cdots X_n^{d_n}$ and $\partial_1^{d_1} \cdots \partial_n^{d_n}$, respectively. 
The symbol $\R_n \coloneqq \bm{R}(X)\langle \partial \rangle$ denotes the noncommutative ring of differential operators with coefficients in $\bm{R}(X)$. 
The subscript $n$ of $\R_n$ is omitted if it is clear from the context. 
We denote the action of an element $l \in \R_n$ on a sufficiently smooth function $\alpha(X) = \alpha(X_1,\dots, X_n)$ as $l \bullet \alpha(X)$; for instance, $\partial_i \bullet \alpha(X) = \pldev{\alpha}{X_i}(X)$. 
We say that a differential operator $l \in \R$ \emph{annihilates} a function $\alpha$ if $l \bullet \alpha = 0$. 
We denote the set of all positive definite $n \times n$ matrices by $\mathrm{PD}(n)$. 
The PDF of a stochastic variable $X$ is denoted by $p(X;Y)$ if it depends on a parameter $Y$. 

\section{Problem Setting \label{sec:problem}}
In the Bayesian filtering approach for discrete-time nonlinear systems, the posterior PDF of the state $x_k \in \bm{R}^n$ at time step $k$ conditional on the output history $y_{[0:k]} \subset \bm{R}^r$ from time steps $0$ to $k$ is recursively updated via the following equations~\cite{Bryson1975}. 
\begin{equation}
    \begin{aligned}
    &p(x_k \mid y_{[0:k]}) = \frac{p(x_k, y_k \mid y_{[0:k-1]})}{p(y_k \mid y_{[0:k-1]})}, \\
    &p_{xy}(x_k, y_k \mid y_{[0:k-1]}) = p(y_k \mid x_k) \\
    & \qquad\quad \times \int_{\bm{R}^n} p(x_k \mid x_{k-1}) p(x_{k-1} \mid y_{[0:k-1]}) dx_{k-1}, \\
    &p(y_k \mid y_{[0:k-1]}) = \int_{\bm{R}^n} p(x_k, y_k \mid y_{[0:k-1]}) dx_k.
    \end{aligned} \label{eq:Bayes}
\end{equation}
We assume that the system dynamics and observation process are defined by the following state and observation equations.  
\begin{align}
    &x_k = f(x_{k-1}, u_k) + w_k, \label{eq:stEq} \\
    &y_k = h(x_k) + v_k, \label{eq:obEq}
\end{align}
where $w_k \in \bm{R}^n$ and $v_k \in \bm{R}^r$ denote the system and observation noises that are assumed to be independent and identically distributed with the PDFs $p_w(w_k)$ and $p_v(v_k)$, respectively. 
In addition, $u_k \in \bm{R}^m$ denotes a given input of the system. 
Using the rule of transformation of PDFs~\cite{Jazwinski1970}, the conditional PDFs $p(x_k \mid x_{k-1})$ and $p(y_k \mid x_k)$ in~\eqref{eq:Bayes} can be rewritten using functions $f$, $h$, $p_w$, and $p_v$ as follows: 
\begin{align}
    p(x_k \mid x_{k-1} ; u_k) &= p_w(x_k - f(x_{k-1}, u_k)), \label{eq:trPDF} \\
    p(y_k \mid x_k) &= p_v(y_k - h(x_k)), \label{eq:obPDF} 
\end{align}
where the former PDF is conditional on $x_{k-1}$ and parametrized by $u_k$ owing to~\eqref{eq:stEq}. 
By substituting~\eqref{eq:trPDF} and~\eqref{eq:obPDF} into~\eqref{eq:Bayes}, we can describe the update law of Bayesian filtering using $f$, $h$, $p_w$, and $p_v$. 

Since Bayesian filtering is a recursive algorithm, we can focus on the one-step update of the posterior PDF and omit the subscript $k$. 
Moreover, the past history of outputs $y_{[0:k-1]}$ has nothing to do with the update at time step $k$; it only appears in the previous posterior PDF $p(x_{k-1} \mid y_{[0:k-1]})$, so we also omit the past history $y_{[0:k-1]}$. 
Thus, \eqref{eq:Bayes} can be summarized as
\begin{equation}
    p(x \mid y; u) = \frac{p_{xy}(x, y ; u)}{\int_{\bm{R}^n} p_{xy}(x, y ; u)dx} , \label{eq:pdf_update}
\end{equation}
where $x$, $y$, and $u$ denote the current state, output, and input, respectively; $x^-$ denotes the previous state, and $p_{xy}(x, y; u)$ denotes the joint PDF of $x$ and $y$ defined as the second line in~\eqref{eq:Bayes}. 

The update law~\eqref{eq:pdf_update} can be viewed as a \emph{functional} that maps the previous posterior PDF $p(x^-)$ to the current one ${p(x \mid y; u)}$. 
In the linear case with Gaussian noise, this functional can be reduced to simple arithmetic operations of the mean and variance of the previous posterior PDF to yield the prediction and update steps of the Kalman filter. 
However, in the nonlinear cases, the Gaussian property of the posterior PDF can no longer be guaranteed, which makes it extremely difficult to compute the functional. 
Therefore, as the first step to overcome this difficulty, we approximate the posterior PDF using a Gaussian. 

For the following discussion, $\tilde{p}$ denotes an approximation of the PDF $p$ with the same stochastic variables. 
Suppose we have a Gaussian $\tilde{p}(x^-; \mu^-, \Sigma^-) = \mathcal{N}(x^-; \mu^-, \Sigma^-)$ approximating $p(x^-)$, where $\mathcal{N}$ is defined as
\[
    \mathcal{N}(\xi; m, V) \coloneqq \frac{1}{\sqrt{(2\pi)^n|V|}}\exp\left( -\frac{1}{2}(\xi-m)^\T V^{-1} (\xi-m) \right).
\]
By replacing $p(x^-)$ with $\tilde{p}(x^- ; \mu^-, \Sigma^-)$ in the definition of $p_{xy}(x, y ; u)$, it can be approximated as 
\begin{multline}
    p_{xy}(x,y ; u) \approx \tilde{p}_{xy}(x, y; u, \mu^-, \Sigma^-) \\\coloneqq p(y \mid x)\int_{\bm{R}^n} p(x \mid x^-; u) \tilde{p}(x^- ; \mu^-, \Sigma^-) dx^-, \label{eq:joint}
\end{multline}
and $p(x \mid y; u)$ can be approximated as  
\begin{equation}
    \tilde{p}(x \mid y; u, \mu^-, \Sigma^-)  = \frac{\tilde{p}_{xy}(x, y ; u, \mu^-, \Sigma^-)}{\int_{\bm{R}^n} \tilde{p}_{xy}(x, y ; u, \mu^-, \Sigma^-) dx}. \label{eq:postPDF}
\end{equation}
Note that the approximating PDF on the left-hand side is \emph{not} Gaussian since we consider the nonlinearities of functions $f$ and $h$ without any approximations. 
Hence, we further approximate $\tilde{p}(x \mid y; u, \mu^-, \Sigma^-)$ with a Gaussian $\mathcal{N}(x ; \mu, \Sigma)$ of the same mean $\mu$ and variance $\Sigma$ as those of $\tilde{p}(x \mid y; u, \mu^-, \Sigma^-)$, that is, $\tilde{p}(x \mid y; u, \mu^-, \Sigma^-) \approx \mathcal{N}(x ; \mu, \Sigma)$ 
where 
\begin{align}
    \mu \coloneqq& \notag E_{\tilde{p}(x \mid y; u, \mu^-, \Sigma^-)}[x] \\
    =& \frac{\int_{\bm{R}^n} x \cdot \tilde{p}_{xy}(x, y ; u, \mu^-, \Sigma^-) dx}{\int_{\bm{R}^n} \tilde{p}_{xy}(x, y ; u, \mu^-, \Sigma^-) dx} \label{eq:mean}
\end{align}
and 
\begin{align}
    \Sigma \coloneqq& E_{\tilde{p}(x \mid y; u, \mu^-, \Sigma^-)}[xx^\T] - \mu\mu^\T \notag \\
    =& \frac{\int_{\bm{R}^n} xx^\T \cdot \tilde{p}_{xy}(x, y ; u, \mu^-, \Sigma^-) dx}{\int_{\bm{R}^n} \tilde{p}_{xy}(x, y ; u, \mu^-, \Sigma^-) dx} - \mu\mu^\T. \label{eq:variance}
\end{align}

Definitions~\eqref{eq:mean} and~\eqref{eq:variance} provide a finite number of functions that map the parameters $y$, $u$, $\mu^-$, and $\Sigma^-$ to the current estimates $\mu$ and $\Sigma$. 
Hence, by evaluating the right-hand sides of~\eqref{eq:mean} and~\eqref{eq:variance} recursively, we can perform Bayesian filtering while exactly considering the nonlinearities of $f$ and $h$. 
For simplicity, we define $\Phi\left[ h(x) \right](y, u, \mu^-, \Sigma^-)$ for a scalar-, vector-, or matrix-valued function $h(x)$ as
\[
    \Phi\left[ h(x) \right](y, u, \mu^-, \Sigma^-) \coloneqq \int_{\bm{R}^n} h(x) \cdot \tilde{p}_{xy}(x, y ; u, \mu^-, \Sigma^-) dx,
\]
which is an integral over $x$ depending on the parameters $y$, $u$, $\mu^-$, and $\Sigma^-$. 
The computations of~\eqref{eq:mean} and~\eqref{eq:variance} are then reduced to evaluations of three integrals, namely $\Phi[1]$, $\Phi[x]$, and $\Phi[xx^\T]$, 
because $\mu$ and $\Sigma$ can be written as $\mu = \Phi[x] / \Phi[1]$ and $\Sigma = \Phi[xx^\T] / \Phi[1] - \mu\mu^\T$, respectively. 

For nonlinear KFs, the posterior PDF is usually approximated as Gaussian~\cite{Arasaratnam2009,Fang2018}.
The integrals $\Phi[1]$, $\Phi[x]$, and $\Phi[xx^\T]$ are also approximated without considering the effects of higher order moments induced by the nonlinearity of $f$ and $h$.  
More specifically, the nonlinear KFs first approximate the integral in~\eqref{eq:joint} by computing the mean and variance of the nonlinear function $f(x^-, u)$ and further approximate the integral in the denominator of~\eqref{eq:postPDF} in the same way using $h(x)$. 
By contrast, the proposed method considers the exact integrals using the HGM. 

On the other hand, the PF approximates the integrals appearing in~\eqref{eq:Bayes} using the Monte-Carlo scheme instead of the Gaussian approximation. 
Therefore, for a sufficiently large number of particles, the PF must exhibit a better accuracy than the proposed method. 
However, in the numerical example provided in Section~\ref{sec:examples}, we observe that the computational cost for the PF is more than that of the proposed method for the same level of accuracy. 

\begin{rem}
    Although, as mentioned above, the Gaussian approximation of the posterior distribution is a common approach in nonlinear filtering theory, this assumption is not suited to the cases of noise with heavy-tailed distributions, such as the Cauchy noise. 
    In such cases, however, the proposed method may still be applicable by approximating the posterior distribution using a heavy-tailed distribution with some parameters. 
    Provided that the parameters of the approximating distribution can be expressed as functions of the previous estimates, observed output, and known input (which corresponds to~\eqref{eq:mean} and~\eqref{eq:variance} for the Gaussian approximation), the proposed method can be performed in almost the same manner, which could be considered in future work. 
\end{rem}


\section{Holonomic Gradient Method}
As described in the previous section, the update law of the mean and variance can be formulated as the evaluation of three integrals $\Phi[1]$, $\Phi[x]$, and $\Phi[xx^\T]$. 
For linear systems with Gaussian noise, these integrals can be computed analytically using the Gaussian integral. 
However, for nonlinear systems, we need a numerical method because it is an extremely complex process to compute the integrals analytically. 
To achieve the trade-off between the exact consideration of nonlinearity and reduction of computational cost, we use a symbolic-numeric method for evaluating the complex integrals, called the HGM. 
We only provide the outline of the HGM in a specific situation due to limited space. 

Let us consider the evaluation of the following integral. 
\begin{equation}
    \alpha(X) \coloneqq \int_{\bm{R}^m} \beta(X, Y) \cdot \gamma(X, Y) dY, \label{eq:comp_int}
\end{equation}
where $\beta$ and $\gamma$ are smooth scalar-valued functions in $X \in \bm{R}^n$ and $Y \in \bm{R}^m$. 
For example, $\Phi[h(x)](y, u, \mu^-, \Sigma^-)$ can be regarded as integral~\eqref{eq:comp_int} with $\beta = h$, $\gamma = \tilde{p}_{xy}$, $Y = x$, and $X$ corresponds to $(y, u, \mu^-, \Sigma^-)$. 
\begin{defi}\label{def:HF}
    The function $\alpha(X)$ is called \emph{holonomic} if there exists a set of differential operators $\mathcal{M} = \{\partial^{d_1},\dots,\allowbreak\partial^{d_{q-1}}\}\ (q \in \bm{Z}_{\geq 0},\ d_i \in \bm{Z}_{\geq 0}^n)$ such that the vector-valued function $Q(X)$ defined as 
    \begin{equation}
        Q = [ \alpha\ \partial^{d_1}\alpha\ \cdots\ \partial^{d_{q-1}} \alpha]^\T
        \label{eq:q}
    \end{equation}
    satisfies a set of the following PDEs. 
    \begin{equation}
        \partial_{X_i} Q(X) = A_{X_i}(X)Q(X) \qquad (i = 1,\dots,n), \label{eq:pf}
    \end{equation}
    where $A_{X_i}(X) \in \bm{R}(X)^{n \times n}$. 
    The set of PDEs~\eqref{eq:pf} is called the \emph{Pfaffian system} of $Q(X)$. 
\end{defi}

We assume that we know $\mathcal{M}$ and $A_{X_i}$ in Definition~\ref{def:HF} for the holonomic function $\alpha(X)$. 
Moreover, we assume that the vector $Q(X_{\mathrm{init}})$ at a certain point $X_{\mathrm{init}}$ is known. 
For example, the vector $Q(X_{\mathrm{init}})$ can be directly computed from~\eqref{eq:comp_int} if we use sufficient computational resources. 
Once $Q(X_{\mathrm{init}})$ is obtained, by virtue of~\eqref{eq:pf}, we can efficiently evaluate $\alpha(\hat{X})$ at a point $\hat{X} \neq X_{\mathrm{init}}$ without directly evaluating the multiple integral~\eqref{eq:comp_int}. 
The target point $\hat{X}$ can be arbitrarily selected considering that the following process is successfully performed. 

Let $X(s)$ be a smooth curve $[0,1] \ni s \mapsto X(s) \in \bm{R}^n$ with $X(0) = X_{\mathrm{init}}$ and $X(1) = \hat{X}$.
Considering the Pfaffian system~\eqref{eq:pf} on this curve, we obtain an IVP
\begin{multline}
   \frac{dQ}{ds} = \sum_{i=1}^n \partial_{X_i}Q \frac{dXi}{ds}= \sum_{i=1}^{n} A_{X_i}(X(s))Q(X(s)) \frac{dX_i}{ds},\\ Q(X(0)) = Q(X_{\mathrm{init}}). \label{eq:ODE}
\end{multline}
This IVP can be solved to obtain $Q(\hat{X})$ using numerical solution methods such as the RK4 method if none of the denominators in $A_{X_i}$ vanish at any $s \in [0,1]$. 
Consequently, we can evaluate $\alpha(\hat{X})$ as the first component of $Q(\hat{X})$. 

The finite set $\mathcal{M}$ and matrices $A_{X_i}$, which are supposed to be known, can be computed from a finite set of differential operators called a \emph{basis of a zero-dimensional ideal} in $\mathcal{R}_n$. 
The following theorem provides a way to find such a basis for a given function and to determine whether the function is holonomic or not (We have skipped the details for simplicity). 
\begin{theo}\cite{Hibi2013} \label{thm:ZD}
    A smooth function $\alpha(X)$ is holonomic if and only if  there exist differential operators, described as finite sums of the following form.
    \begin{equation}
        l_i \coloneqq \sum_{d} a_d(X)\partial_{X_i}^d \quad (a_d \in \bm{R}(X),\ i=1,\dots,n), \label{eq:ZDbasis}
    \end{equation}
    which annihilate $\alpha(X)$, i.e., $l_i \bullet \alpha(X) = 0$ for $(i = 1,\dots,n)$. 
    Moreover, the set of differential operators $\{l_1,\dots,l_n\}$ is a basis of a zero-dimensional ideal in $\mathcal{R}_n$. 
\end{theo}
\begin{ex}
    A function $\alpha(X_1, X_2) = \cos(X_1X_2)$ is holonomic because there exist two differential operators of the form~\eqref{eq:ZDbasis}: 
    \[
        l_1 \coloneqq \partial_{X_1}^2 + X_2^2, \qquad l_2 \coloneqq \partial_{X_2}^2 + X_1^2, 
    \]
    which annihilate $\alpha(X_1, X_2)$. Indeed, by using, for example, \texttt{Risa/Asir}~\cite{NoroAsir} and its library \texttt{yang}, the Pfaffian system~\eqref{eq:pf} can be computed as 
    \[
        \partial_{X_1}Q = \left[
            \begin{array}{cc}
                1/X_1 & X_1X_2 \\
                X_2/X_1 & 0
            \end{array}
        \right] Q, \quad 
        \partial_{X_2}Q = \left[
            \begin{array}{cc}
                0 & X_1^2 \\
                1 & 0
            \end{array}
        \right] Q,  
    \] where $Q(X_1, X_2) = [\alpha\ \partial_{X_2}\alpha]^\T$, that is, $\mathcal{M}$ in Definition~\ref{def:HF} is a singleton $\{\partial_{X_2}\}$. 
\end{ex}

It may be difficult to find a basis of a zero-dimensional ideal for $\alpha$ directly from the expression~\eqref{eq:comp_int}. 
Using symbolic computation, we can compute a basis of a zero-dimensional ideal for $\alpha$ from those for $\beta$ and  $\gamma$~\cite{Oaku2003}, which implies the following corollary~\cite{Oaku2003,Hibi2013}. 
\begin{cor}\label{cor:CP}
    If $\beta$ and $\gamma$ in~\eqref{eq:comp_int} are holonomic, 
    the product $\beta \cdot \gamma$ is also holonomic. 
    In addition, if the product $\beta \cdot \gamma$ is rapidly decreasing with respect to each $Y_i\ (i = 1,\dots,m)$, 
    then $\alpha$, defined as~\eqref{eq:comp_int}, is also holonomic. 
\end{cor}
\begin{ex}
    For $X \in \bm{R}$ and $Y \in \bm{R}$, let $\beta(X, Y) = \cos(XY)$ and $\gamma(X, Y) = \exp(-X^2/2 - Y^2/2)$. 
    Evidently, $\beta$ is annihilated by $\partial_X^2 + Y^2$ and $\partial_Y^2 + X^2$, 
    and $\gamma$ is annihilated by $\partial_X + X$ and $\partial_Y + Y$. 
    That is, both functions are holonomic. 
    Moreover, $\gamma$ is rapidly decreasing with respect to $Y$ and so is the product $\beta \cdot \gamma$, which indicates that $\alpha(X)$ defined by~\eqref{eq:comp_int} is holonomic (Corollary~\ref{cor:CP}). 
    Using \texttt{Risa/Asir} and its library \texttt{nk\_restriction}, we can compute $\partial_X + 2X$ as a basis of a zero-dimensional ideal for $\alpha(X)$ from the differential operators annihilating $\beta$ and $\gamma$. 
    Indeed, it is easy to confirm that in this example $\alpha(X) \propto \exp(-X^2)$, which is annihilated by the basis.
\end{ex}

Corollary~\ref{cor:CP}, combined with the following assumption, ensures that $\tilde{p}_{xy}(x, y ; u, \mu^-, \Sigma^-)$ is holonomic. 
\begin{assm} \label{assm:HF}
    The conditional PDFs $p(x \mid x^-; u)$ and $p(y \mid x)$ in~\eqref{eq:joint} are holonomic functions. 
\end{assm}

\section{Evaluation of Mean and Variance via Integral Transform}
Under Assumption~\ref{assm:HF}, all the components of the integrals $\Phi[1]$, $\Phi[x]$, and $\Phi[xx^\T]$ are holonomic functions of $u$, $y$, $\mu^-$, and $\Sigma^-$.
Therefore, they can be efficiently evaluated using the HGM if the corresponding Pfaffian systems~\eqref{eq:pf} and initial vectors~\eqref{eq:q} are computed in advance. 
The authors propose a method to compute~\eqref{eq:mean} and~\eqref{eq:variance} by directly using the HGM~\cite{Iori2020}. 
Although the previous method showed better performance than other existing methods for a numerical example, the IVP~\eqref{eq:ODE} has to be solved for each component of the three integrals. 
If we can decrease the number of IVPs, then it reduces the online computational cost. 
Hereafter, $z \in \bm{R}^N$ denotes a vector consisting of all the independent components of $y \in \bm{R}^r$, $u \in \bm{R}^m$, $\mu^- \in \bm{R}^n$, and $\Sigma^- \in \mathrm{PD}(n)$, where $N = r + m + n + n(n+1)/2$. 
With this notation, for example, $\tilde{p}_{xy}(x, z)$ denotes $\tilde{p}_{xy}(x, y ; u, \mu^-, \Sigma^-)$ as a function of $x$, $y$, $u$, $\mu^-$, and $\Sigma^-$. 

In this work, we use an integral transform of $\tilde{p}_{xy}$ to reduce the number of IVPs solved in the one-step estimation. 
Consider an integral transform $\mathcal{T}[\tilde{p}_{xy}]$ defined as follows.
\begin{equation}
    \mathcal{T}[\tilde{p}_{xy}](\xi, z) \coloneqq \int_{\bm{R}^n} \exp(\xi^\T x) \cdot \tilde{p}_{xy}(x, z)dx. \label{eq:IT}
\end{equation}
Note that although this definition is similar to that of the moment generating function, it is different because $\tilde{p}_{xy}$ is the PDF of the $x$ and $y$ pair rather than $x$. 
We assume that there exists a compact subset of $\bm{R}^n$ including the origin such that for all $\xi$ in the subset, \eqref{eq:IT} can be defined and smooth. 

The integrals $\Phi[1]$, $\Phi[x]$, and $\Phi[xx^\T]$ are then obtained from the derivatives of $\mathcal{T}[\tilde{p}_{xy}]$ at $\xi = 0$ as follows: 
\begin{equation}
    \mathcal{T}[\tilde{p}_{xy}](0, z) \notag = \int_{\bm{R}^n} \tilde{p}_{xy}(x, z)dx = \Phi[1](z), 
\end{equation}
\begin{align}
    \partial_{\xi_i} \mathcal{T}[\tilde{p}_{xy}](0, z) &= \int_{\bm{R}^n} \left.\left[\partial_{\xi_i} \exp(\xi^\T x)\tilde{p}_{xy}(x, z)\right]\right|_{\xi=0}dx \notag \\
    &= \Phi\left[x_i\right](z), \label{eq:ITd1}
\end{align}
\begin{align}
    \partial_{\xi_i}\partial_{\xi_j}\mathcal{T}[\tilde{p}_{xy}](0, z) &= \int_{\bm{R}^n} \left.\left[\partial_{\xi_i}\partial_{\xi_j} \exp(\xi^\T x)\tilde{p}_{xy}(x, z)\right]\right|_{\xi=0}dx \notag \\
    &= \Phi\left[x_ix_j\right](z). \label{eq:ITd2}
\end{align}

Here, the approximating PDF $\tilde{p}_{xy}$ is a holonomic function under Assumption~\ref{assm:HF}, and the kernel $\exp(\xi^\T x)$ is also holonomic. 
This implies that $\mathcal{T}[\tilde{p}_{xy}]$ is also a holonomic function by Corollary~\ref{cor:CP}. 
Therefore, $\mathcal{M}$ exists such that the vector-valued function $Q(X)$ defined by $\mathcal{M}$ and $\mathcal{T}[\tilde{p}_{xy}]$ satisfies the Pfaffian system~\eqref{eq:pf}. 
Using the Pfaffian system, we can explicitly express $\mathcal{T}[\tilde{p}_{xy}]$ and its derivatives~\eqref{eq:ITd1} as well as~\eqref{eq:ITd2} as a linear combination of the components of $Q$ with coefficients in $\bm{R}(\xi, z)$. 

First, $\mathcal{T}[\tilde{p}_{xy}](\xi, z)$ itself is the first component of $Q(\xi, z)$, thus satisfying the following: 
\begin{equation}
    \Phi[1] = \mathcal{T}[\tilde{p}_{xy}] = C^{(0)}Q, \label{eq:C0}
\end{equation}
where $C^{(0)}$ is a constant coefficient vector $[1\ 0\ \cdots\ 0] \in \bm{R}^q$. 
Next, the first derivatives $\partial_{\xi_i} \mathcal{T}[\tilde{p}_{xy}]\ (i=1,\dots,n)$ are obtained using the first components of the left-hand sides of~\eqref{eq:pf}. 
Hence, the following identities hold. 
\begin{equation}
    \Phi[x_i] = \partial_{\xi_i} \mathcal{T}[\tilde{p}_{xy}] =  C^{(1)}_iQ \quad (i=1,\dots,n), \label{eq:C1}
\end{equation}
where the coefficient vector $C^{(1)}_i(\xi, z)$ is the first row vector of $A_{\xi_i}(\xi, z) \in \bm{R}(\xi, z)^{q \times q}$ in~\eqref{eq:pf}. 
Finally, the second derivatives $\partial_{\xi_i} \partial_{\xi_j} \mathcal{T}[\tilde{p}_{xy}]\ (i, j=1,\dots,n)$ are obtained by differentiating both sides of~\eqref{eq:pf}; $\partial_{\xi_i} \bullet \partial_{\xi_j} Q = \partial_{\xi_i} \bullet (A_{\xi_j}Q) = \partial_{\xi_i}A_{\xi_j}Q + A_{\xi_j}\partial_{\xi_i}Q = (\partial_{\xi_i}A_{\xi_j} + A_{\xi_j}A_{\xi_i})Q$, 
and hence 
\begin{equation}
    \Phi[x_i x_j] = \partial_{\xi_i} \partial_{\xi_j} \mathcal{T}[\tilde{p}_{xy}] = C^{(2)}_{ij}Q \quad (i, j \in \{1,\dots,n\}), \label{eq:C2}
\end{equation}
where the coefficient vector $C^{(2)}_{ij}(\xi, z)$ is the first row vector of $\partial_{\xi_j}A_{\xi_i} + A_{\xi_i}A_{\xi_j}$. 

We assume the explicit expressions of $C^{(1)}_i(\xi, z)$ and $C^{(2)}_{ij}(\xi, z)$ are obtained. 
Using the expressions, we can compute the values of $\Phi[1]$, $\Phi[x]$, and $\Phi[xx^\T]$ at $z = \hat{z}$ from~\eqref{eq:C0}--\eqref{eq:C2} if we have the vector $Q(0, \hat{z})$ for a given $\hat{z}$. 
The computation of $Q(0, \hat{z})$ can be accomplished by numerically integrating~\eqref{eq:ODE} using numerical integration methods, such as the RK4 method. 
Therefore, the current mean $\mu$ and variance $\Sigma$ for a given data $\hat{z}$ can be computed from~\eqref{eq:mean} and~\eqref{eq:variance} as 
\begin{equation}
    \mu = \frac{\Phi[x](\hat{z})}{\Phi[1](\hat{z})}, \quad \Sigma = \frac{\Phi[xx^\T](\hat{z})}{\Phi[1](\hat{z})} - \mu\mu^\T. \label{eq:estimates}
\end{equation}
In the previous method~\cite{Iori2020}, IVP~\eqref{eq:ODE} has to be solved for each of the three integrals $\Phi[1]$, $\Phi[x]$, and $\Phi[xx^\T]$ even for the one-dimensional case, that is, three IVPs have to be solved in total. 
By contrast, only a single IVP for $\mathcal{T}[\tilde{p}_{xy}]$ needs to be solved in the proposed method. 
This decrease in the number of IVPs leads to a reduction in the computational cost, which can be seen in Section~\ref{sec:examples}. 
\begin{algorithm}[b]
    \caption{Computation of Pfaffian System for $\mathcal{T}[\tilde{p}_{xy}]$ (Offline part)\label{alg:offline}} 
    \begin{algorithmic}[1]
        \Require{Bases of zero-dimensional ideals annihilating $\tilde{p}(x^- ; \mu^-, \Sigma^-)$, $p(x \mid x^-; u)$, $p(y \mid x)$, and $\exp(\xi^\T x)$}
        \Ensure{Explicit expressions of $C^{(1)}_i(\xi, z)\ (i=1,\dots,n)$, $C^{(2)}_{ij}(\xi, z)\ (i, j \in \{1,\dots,n\})$, $A_{\lambda}(\xi, z)\ (\lambda \in \{\xi_1,\dots,\xi_n,z_1,\dots,z_N\})$}
        \State{Compute basis $\mathcal{B}$ of zero-dimensional ideal for $\mathcal{T}[\tilde{p}_{xy}]$ from inputs \label{algline:computeB}}
        \State{Find $\mathcal{M} = \{\partial^{d_1},\dots,\partial^{d_{q-1}}\}$ in Definition~\ref{def:HF} from $\mathcal{B}$ \label{algline:computeM}}
        \State{Compute coefficient matrices $A_\lambda$ from $\mathcal{B}$ and $\mathcal{M}$ \label{algline:computeA}}
        \State{Compute $C^{(1)}_i(\xi, z)$ and $C^{(2)}_{ij}(\xi, z)$ from $A_{\xi_i}$ \label{algline:computeC}}
    \end{algorithmic}
\end{algorithm}
\begin{algorithm}[b]
    \caption{One-step Estimation via HGM using~\eqref{eq:IT} (Online part)\label{alg:SE_HGM_IT}}
    \begin{algorithmic}[1]
        \Require{Explicit expressions of $C^{(0)}$, $C^{(1)}_i(\xi, z)\ (i=1,\dots,n)$, $C^{(2)}_{ij}(\xi, z)\ (i, j \in \{1,\dots,n\})$, $A_{\lambda}(\xi, z)\ (\lambda \in \{\xi_1,\dots,\xi_n,z_1,\dots,z_N\})$, initial point $z_{\mathrm{init}}$ and vector $Q(0, z_{\mathrm{init}})$, as well as given data $\hat{z}$}
        \Ensure{Failure or estimates $\mu$ and $\Sigma$}
        \State{Define $z(s)$ such that $z(0) = z_{\mathrm{init}}$ and $z(1) = \hat{z}$ \label{algline:path}}
        \If{$\exists s \in [0,1]$ at which denominator in $A_\lambda$ vanishes,}
            \State \Return{algorithm has failed \label{algline:fail}}
        \EndIf
        \State{Solve IVP~\eqref{eq:ODE} numerically using initial vector $Q(0, z_{\mathrm{init}})$ and obtain $Q(0, \hat{z})$ \label{algline:integrate}}
        \State \Return{\eqref{eq:estimates} obtained from~\eqref{eq:C0}--\eqref{eq:C2}}
    \end{algorithmic}
\end{algorithm}

Now, we summarize the offline and online processes of the proposed method as Algorithms~\ref{alg:offline} and~\ref{alg:SE_HGM_IT}. 
The details of each step in Algorithm~\ref{alg:offline} can be found in~\cite{Oaku2003,Iori2020}. 

\section{Numerical Example \label{sec:examples}}
This section presents a numerical example to show the efficiency of the proposed method. 
In the following demonstration, we use \texttt{Risa/Asir}~\cite{NoroAsir} for the symbolic computation, \texttt{Maple} for the offline numerical computation, and \texttt{Python} for the online numerical computation on a PC (Intel(R) Core(TM) i7-1065G7 CPU @ 1.30 GHz; RAM: 16 GB). 

\subsection{Problem setting and results of offline computation}
Consider the following nonlinear stochastic one-dimensional system. 
\begin{equation}
    x =  \frac{4}{5}x^- + u + w, \quad y = \frac{2x}{1 + x^2} + v, \label{eq:exsys}
\end{equation}
where $w$ and $v$ are the system and observation noises with the standard Gaussian PDF, respectively, and $u$ is given as the function $\cos(0.6k)$ of the discrete-time step $k$. 

First, we symbolically compute the inputs of Algorithm~\ref{alg:offline}, namely, bases of zero-dimensional ideals annihilating ${\tilde{p}(x^- ; \mu^-, \Sigma^-)} = \mathcal{N}(x^-; \mu^-, \Sigma^-)$, \eqref{eq:trPDF}, \eqref{eq:obPDF}, and $\exp{(\xi x)}$. 
By differentiating $\tilde{p}(x^- ; \mu^-, \Sigma^-)$ with respect to each variable, we can obtain the differential operators annihilating $\tilde{p}(x^- ; \mu^-, \Sigma^-)$ as
\begin{equation}
    \begin{aligned}
        & \Sigma^-\partial_{x^-} + x^- - \mu^-, \quad \Sigma^-\partial_{\mu^-} - x^- + \mu^-, \\
        & 2(\Sigma^-)^2\partial_{\Sigma^-} + \Sigma^- - (x^- - \mu^-)^2. 
    \end{aligned} \label{eq:genIn}
\end{equation}
The conditional PDF~\eqref{eq:trPDF} is obtained by substituting the state equation into the PDF of $\mathcal{N}(w;0, 1)$ and is annihilated by the following three differential operators:
\begin{equation}
    \begin{aligned}
    & \partial_x + x - (4/5)x^- - u, \quad \partial_u -x + (4/5)x^- + u, \\
    & -(5/4)\partial_{x^-} + x - (4/5)x^- - u. 
    \end{aligned} \label{eq:genTr}
\end{equation}
The other conditional PDF~\eqref{eq:obPDF} is also obtained by substituting the observation equation into the PDF of $\mathcal{N}(v;0, 1)$ and is annihilated by 
\begin{equation}
    \partial_y + y - \frac{2x}{1 + x^2}, \quad -\frac{(1 + x^2)^2}{2(1-x^2)}\partial_x + y - \frac{2x}{1 + x^2}. \label{eq:genOb}
\end{equation}
Finally, $\exp(\xi x)$ is annihilated by 
\begin{equation}
    \partial_\xi - x, \quad \partial_x - \xi. \label{eq:genK}
\end{equation}
By Theorem~\ref{thm:ZD}, the sets of differential operators~\eqref{eq:genIn}--\eqref{eq:genK} are bases of zero-dimensional ideals annihilating $\tilde{p}(x^-; \mu^-, \Sigma^-)$, $p(x \mid x^-; u)$, $p(y \mid x)$, and $\exp(\xi x)$, respectively. 

Using Algorithm~\ref{alg:offline}, we can derive a basis $\mathcal{B}$ of a zero-dimensional ideal annihilating $\mathcal{T}[\tilde{p}_{xy}]$, a finite set of differential operators $\mathcal{M}$, and the coefficient vectors $C^{(1)}_i(\xi, z), C^{(2)}_{ij}(\xi, z)$ and matrices $A_\lambda\ (\lambda \in \{\xi, z_1, z_2, z_3,z_4\})$, where $z = [z_1\ \cdots \ z_4]^\T = [y\ u\ \mu^-\ \Sigma^-]^\T$. 
In this example, $\mathcal{B}$ consists of seven differential operators, which include
\[
    2\partial_\xi\partial_{z_1}^2 - \partial_\xi^2\partial_{z_2} + 2y\partial_\xi\partial_{z_1} + \xi\partial_\xi^2 - 2\partial_{z_1} - \partial_{z_2} + 4\partial_\xi + \xi. 
\]
The other operators are omitted because of space limitations. 
The finite set $\mathcal{M}$ is obtained as $\left\{\partial_{z_1}, \partial_{z_2}, \partial_{z_4}, \partial_{z_1}^2, \partial_{z_2}\partial_{z_1}, \partial_{z_4}\partial_{z_2}\right\}$. 
Hence, $Q$ becomes a seven-dimensional vector-valued function: 
\begin{equation}
    Q = \left[1\ \partial_{z_1}\ \partial_{z_2}\ \partial_{z_4}\ \partial_{z_1}^2\ \partial_{z_2}\partial_{z_1}\ \partial_{z_4}\partial_{z_2}\right]^\T \bullet \mathcal{T}[\tilde{p}_{xy}], \label{eq:exQ}
\end{equation}
where $\bullet$ represents the actions of all the components on $\mathcal{T}[\tilde{p}_{xy}]$. 
All entries of the coefficient matrices $A_\lambda \in \bm{R}(\xi, z)^{7 \times 7}\ (\lambda \in \{\xi, z_1, z_2, z_3, z_4\})$ as well as the coefficient vectors in~\eqref{eq:C1} and~\eqref{eq:C2} are explicitly computed from $\mathcal{B}$. 
For example, an overview of $A_\xi$ is presented bellow. 
\[
    \left[
        \begin{array}{ccccccc}
            z_2 + \frac{4}{5}z_3 & 0 & \frac{16}{25}z_4 + 1 & 0 & 0 & 0 & 0 \\
            0 & z_2 + \frac{4}{5}z_3 & 0 & 0 & 0 & * & 0 \\
            1 & 0 & z_2 + \frac{4}{5}z_3 & * & 0 & 0 & 0 \\
            0 & 0 & \frac{16}{25} & * & 0 & 0 & * \\
            * & (-z_2-\frac{4}{5}z_3)z_1+1 & * & * & 0 & * & * \\
            * & * & * & * & 0 & * & 0 \\
            * & * & * & * & * & * & * 
        \end{array}
    \right],
\]
where the entries $*$ denote rational functions of $\xi,z_1,\dots,z_4$ and are omitted because of the space limitations. 

Finally, we have to choose the initial points $z_{\mathrm{init}}$ and compute corresponding initial vectors $Q(0, z_{\mathrm{init}})$. 
The choice of the initial points depends on the denominators in $A_\lambda\ (\lambda \in \{\xi, z_1, z_2, z_3, z_4\})$ because they should not vanish on the integration path; otherwise, the right-hand side of the ODE in~\eqref{eq:ODE} diverges and therefore Algorithm~\ref{alg:SE_HGM_IT} fails (line:~\ref{algline:fail}). 
In this example, however, the least common multiple of all the denominators in $A_\lambda$ is $16z_4 + 25$, which must not be zero because $z_4 = \Sigma^- > 0$. 
Therefore, any component of $A_\lambda$ must not diverge. 
For this example, we choose the following eight initial points: $Z_{\mathrm{init}} \coloneqq \left\{[z_1\ z_2\ z_3\ z_4]^\T \mid z_1 = \pm 1, z_2 = \pm 1, z_3 = \pm 1, z_4 = 1 \right\}$. 
The corresponding initial vectors can be computed numerically from the definition~\eqref{eq:exQ}. 

\subsection{Settings for online computation and estimation results}
\begin{figure}[t]
    \centering
    \begin{minipage}[b]{\columnwidth}
        \centering
        \includegraphics[width=0.95\columnwidth]{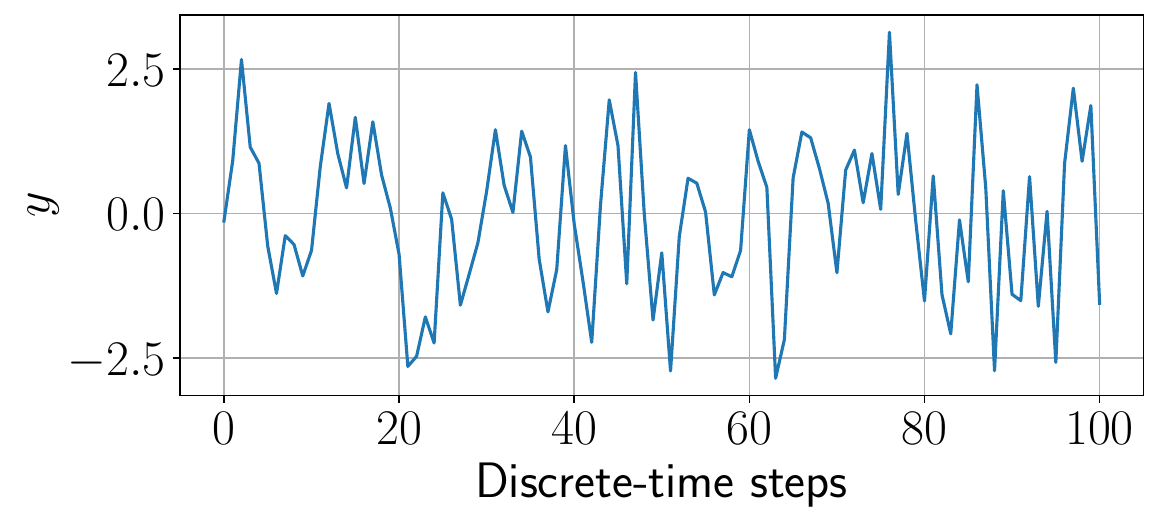}
        \subcaption{Realization of output sequence}
    \end{minipage}
    \begin{minipage}[b]{\columnwidth}
        \centering
        \includegraphics[width=0.95\columnwidth]{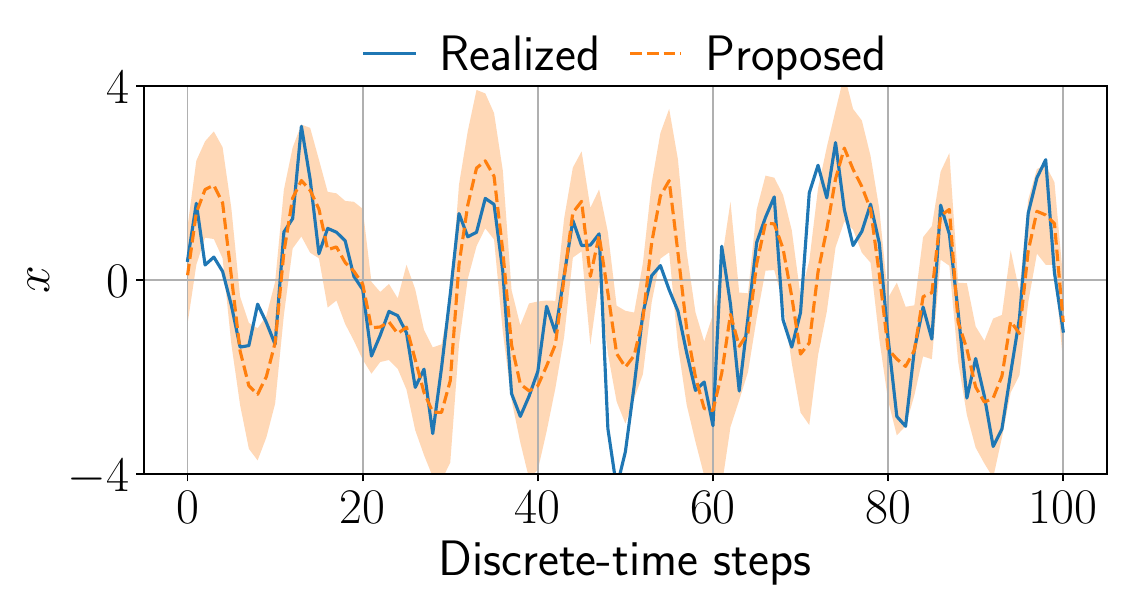}
        \subcaption{Trajectories of realized state (solid) and estimate (dashed). Filled area shows intervals between $\mu \pm \sqrt{\Sigma}$}
    \end{minipage}
    \caption{Realization of system~\eqref{eq:exsys} and corresponding estimation result \label{fig:est}}
\end{figure}

In the online computations, the data $z$ consisting of the input $u$, output $y$, and previous estimates $\mu^-$ and $\Sigma^-$ are given at each time step. 
Let $\hat{z}$ denote the given data at a time step. 
For this example, we fix the integration path in line~\ref{algline:path} of Algorithm~\ref{alg:SE_HGM_IT} to a line segment from an initial point $z_{\mathrm{init}}$ to a given data $\hat{z}$, that is,  $z(s) \coloneqq s\hat{z} + (1 - s)z_{\mathrm{init}}$. 
The initial point $z_{\mathrm{init}}$ is therefore selected for each given data $\hat{z}$ such that the length of the line segment is the minimum among the prescribed set $Z_{\mathrm{init}}$. 
The numerical integration along the defined path is performed using the ABM4 method, with the first three steps initialized using the RK4 method. 
Under these settings, Algorithm~\ref{alg:SE_HGM_IT} is performed at each time step to compute the next estimates $\mu$ and $\Sigma$ from the previous ones $\mu^- (=z_3)$ and $\Sigma^- (=z_4)$ using the measurement $y (=z_1)$ and the known input $u (=z_2)$. 
Figure~\ref{fig:est} shows a realization of the system~\eqref{eq:exsys} and the estimation results obtained from the proposed method. 

For comparison, we implemented the EKF, UKF, PF, and our previous method. 
The number of particles in the PF was set to 100 so that its accuracy would be almost the same as that of the proposed method.  
Three hundred realizations are generated from the random initial states and sequences of system and observation noises for this comparison. 

\begin{figure}[t]
    \centering
    \includegraphics[width=0.95\columnwidth]{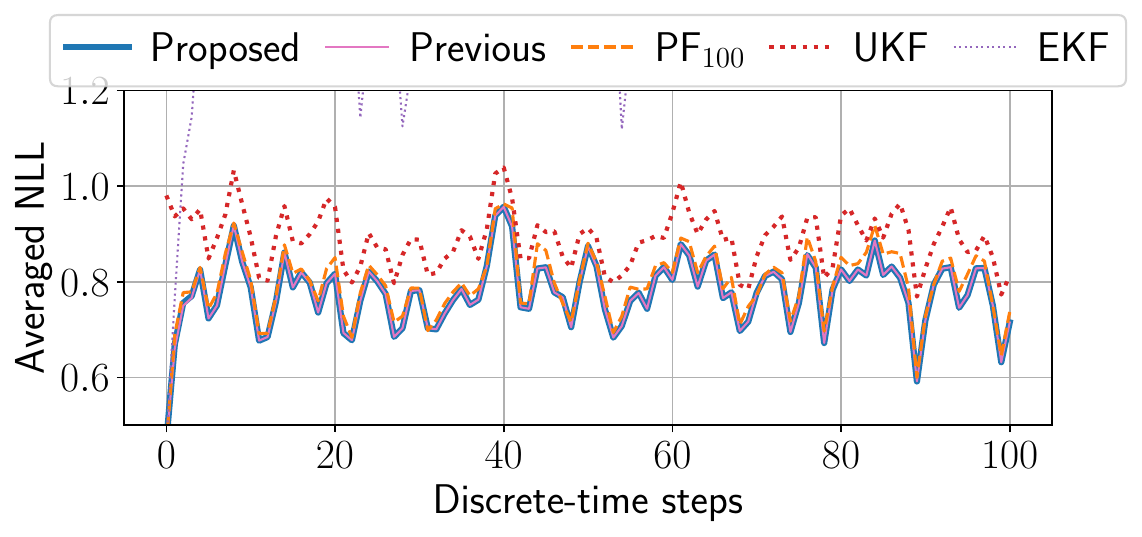}
    \caption{Comparison of NLLs for proposed method (solid, thick), previous method (solid, thin), PFs of 100 particles (dashed), UKF (dotted), and EKF (dash-dotted, almost greater than 1.2 after $k = 4$) \label{fig:nll}}
\end{figure}

\begin{figure}[t]
    \centering
    \includegraphics[width=0.75\columnwidth]{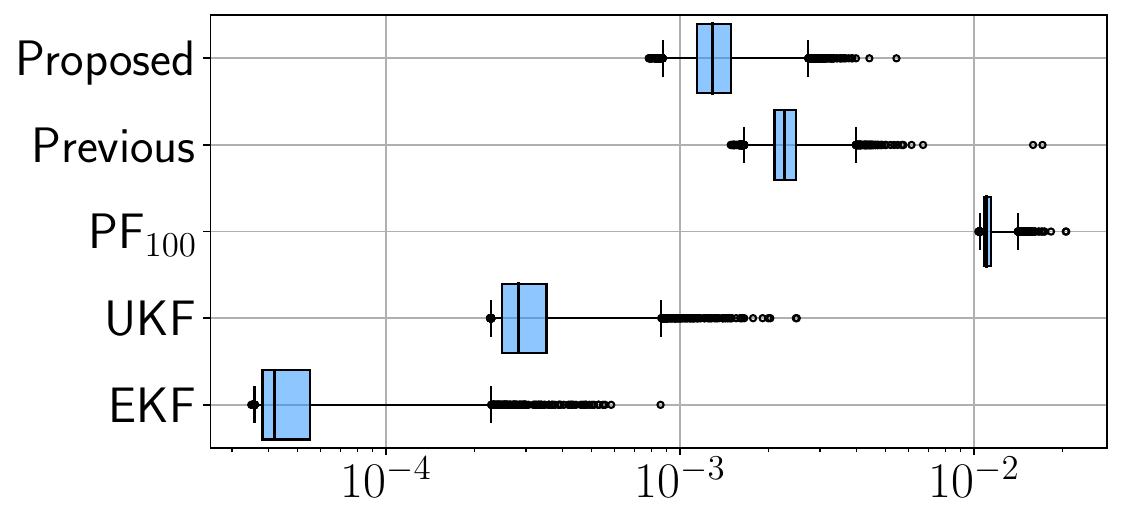}
    \caption{Boxplots of computational times for one-step estimations of all methods \label{fig:cptimes}}
\end{figure}

To compare the performances of all the methods, the negative log-likelihood (NLL)~\cite{Deisenroth2010} was computed.  
Figure~\ref{fig:nll} shows the NLL of each method averaged over all realizations. 
As can be seen, only the PF of 100 particles shows comparable performance to the proposed and previous methods. 
The averaged NLL of the previous method is identical to that of the proposed method because both approaches compute the same estimates~\eqref{eq:estimates}. 
The computational times for the one-step estimations of all the methods are summarized as boxplots in Fig.~\ref{fig:cptimes}. 
It can be seen that the proposed method is faster than the previous method as well as the PF with 100 particles, which show comparable performances for the NLL values. 
These results indicate that the proposed method outperforms the PF and our previous method.

As shown in this example, there exists at least a single problem in which the proposed method outperforms the existing methods, exhibiting the superiority of the proposed method. 
Currently, it is difficult to demonstrate the superiority for general cases because of limitations such as the instability of ODE~\eqref{eq:ODE}. 
Therefore, the theoretical or experimental analysis demonstrating the superiority will be a part of future research. 

Despite the complex procedure, the principle behind the proposed method is simple; it is basically the evaluation of the integrals $\Phi[1]$, $\Phi[x]$, and $\Phi[xx^\T]$, which can be naturally obtained by the Gaussian approximation of the posterior PDF in Bayesian filtering.
Thus, the proposed method is simple to understand except for the precise evaluation of integrals, which is complex. 
Moreover, if it is formulated as a set of Algorithms~\ref{alg:offline} and~\ref{alg:SE_HGM_IT}, the implementation of the proposed method can be automated using symbolic computation so that the users do not have to be concerned about its complex procedure. 
Therefore, we believe that it is worth studying the proposed method, and the improvement of its usability should be a part of the future research. 


\section{Conclusion \label{sec:conclusion}}
A symbolic-numeric method is proposed herein to perform Bayesian filtering for a class of nonlinear stochastic systems. 
The posterior PDF of the state is first approximated using a Gaussian PDF, and the mean and variance are then updated while considering the nonlinearity of the system exactly. 
This update process is formulated as evaluations of several integrals using the HGM. 
We also propose an integral transform to compute the mean and variance efficiently from the vector-valued function obtained using the HGM. 

It is known that the numerical integration in the third step of the HGM may diverge owing to errors in the initial vectors. 
Therefore, in future work, the selection of integration paths or initial points will be studied in depth to enable more stable numerical integrations, for example, by considering the special structure of the problems considered in this work. 
Another direction for future work is approximating the posterior PDF more accurately using a PDF with more parameters than the Gaussian, such as the skew Gaussian. 

%
%


\begin{thebibliography}{xxx}


\bibitem{Kalman1960b}R.~E. Kalman, ``{A new approach to linear filtering and
  prediction problems},'' \textit{ ASME Journal of Basic Engineering}, vol. 82,
  pp. 35--45, 1960.

\bibitem{Kalman1961}R.~E. Kalman and R.~S. Bucy, ``{New results in linear
  filtering and prediction theory},'' \textit{ ASME Journal of Basic
  Engineering}, vol. 83, no. 1, pp. 95--108, 1961.

\bibitem{Rusnak2015}I.~Rusnak, ``{Maximum likelihood optimal estimator of
  non-autonomous nonlinear dynamic systems},'' in  \textit{Proceedings of
  European Control Conference}, pp. 909--914, 2015.

\bibitem{Luo2015}X.~Luo, Y.~Jiao, W.~L. Chiou, and S.~S. Yau, ``{A novel
  suboptimal method for solving polynomial filtering problems},'' \textit{
  Automatica}, vol. 62, pp. 26--31, 2015.

\bibitem{Chen2018}B.~Chen and G.~Hu, ``{Nonlinear state estimation under
  bounded noises},'' \textit{ Automatica}, vol. 98, pp. 159--168, 2018.

\bibitem{Duong2018}N.~Duong, J.~L. Speyer, J.~Yoneyama, and M.~Idan, ``{Laplace
  estimator for linear scalar systems},'' in  \textit{Proceedings of the IEEE
  Conference on Decision and Control}, pp. 2283--2290, 2018.

\bibitem{Li2021}K.~Li, F.~Pfaff, and U.~D. Hanebeck, ``{Unscented dual
  quaternion particle filter for SE(3) estimation},'' \textit{ IEEE Control
  Systems Letters}, vol. 5, no. 2, pp. 647--652, 2021.

\bibitem{Kitagawa1993}G.~Kitagawa, ``{Monte Carlo filtering and smoothing
  method for non-Gaussian nonlinear state space model},'' \textit{ Institute of
  Statistical Mathematics Research Memorandum}, vol. 462, 1993.

\bibitem{Gordon1993}N.~J. Gordon, D.~J. Salmond, and A.~F.~M. Smith, ``{Novel
  approach to nonlinear/non-Gaussian Bayesian state estimation},'' in
  \textit{IEE Proceedings F - Radar and Signal Processing}, vol. 140, pp.
  107--113, 1993.

\bibitem{Jazwinski1970}A.~H. Jazwinski \textit{ {Stochastic Processes and
  Filtering Theory}}, Elsevier Science, 1970.

\bibitem{Julier1997}S.~J. Julier and J.~K. Uhlmann, ``{New extension of the
  Kalman filter to nonlinear systems},'' in  \textit{SPIE}, vol. 3068, 1997.

\bibitem{Arasaratnam2009}I.~Arasaratnam and S.~Haykin, ``{Cubature Kalman
  filters},'' \textit{ IEEE Transactions on Automatic Control}, vol. 54, no. 6,
  pp. 1254--1269, 2009.

\bibitem{Ito2000}K.~Ito and K.~Xiong, ``{Gaussian filters for nonlinear
  filtering problems},'' \textit{ IEEE Transactions on Automatic Control}, vol.
  45, no. 5, pp. 910--927, 2000.

\bibitem{Iori2020}T.~Iori and T.~Ohtsuka, ``{Symbolic-numeric computation of
  posterior mean and variance for a class of discrete-time nonlinear stochastic
  systems},'' in  \textit{Proceedings of the IEEE Conference on Decision and
  Control}, pp. 4814--4821, 2020.

\bibitem{Nakayama2011}H.~Nakayama, K.~Nishiyama, M.~Noro, K.~Ohara, T.~Sei,
  N.~Takayama, and A.~Takemura, ``{Holonomic gradient descent and its
  application to the Fisher-Bingham integral},'' \textit{ Advances in Applied
  Mathematics}, vol. 47, no. 3, pp. 639--658, 2011.

\bibitem{Idan2014}M.~Idan and J.~L. Speyer, ``{Multivariate Cauchy estimator
  with scalar measurement and process noises},'' \textit{ SIAM Journal on
  Control and Optimization}, vol. 52, no. 2, pp. 1108--1141, 2014.

\bibitem{Bryson1975}J.~A.~E. Bryson and Y.-C. Ho \textit{ {Applied Optimal
  Control}}, John Wiley \& Sons, 1st edition, 1975.

\bibitem{Fang2018}H.~Fang, N.~Tian, Y.~Wang, M.~Zhou, and M.~A.~Haile, 
``{Nonlinear Bayesian estimation: from Kalman filter to a broader horizon},'' \textit{IEEE/CAA Journal of Automatica Sinica}, vol. 5, no. 2, pp.~401--417, 2018.

\bibitem{Hibi2013}T.~Hibi ed.  \textit{ {Gr{\"{o}}bner Bases: Statistics and
  Software Systems}}, Springer Japan, 1st edition, 2013.

\bibitem{Oaku2003}T.~Oaku, Y.~Shiraki, and N.~Takayama, ``{Algebraic algorithms
  for D-modules and numerical analysis},'' \textit{ Computer Mathematics
  (Proceedings of ASCM 2003)}, vol. 10, pp. 23--39, 2003.

\bibitem{NoroAsir}M.~Noro, N.~Takayama, H.~Nakayama, K.~Nishiyama, and
  K.~Ohara, ``{Risa/Asir: A computer algebra system},'' \url{http://www.math.kobe-u.ac.jp/Asir/asir.html}, 2020.

\bibitem{Deisenroth2010}M.~P. Deisenroth \textit{ {Efficient Reinforcement
  Learning Using Gaussian Processes}}, KIT Scientific Publishing, 2010.

\end{thebibliography}

\end{document}